\documentclass[12pt]{amsart}
\usepackage{a4wide,bm}
\usepackage{amssymb}
\usepackage{euscript}
\usepackage{xcolor}
\usepackage{geometry}
\usepackage{microtype}
\usepackage{graphicx}
\usepackage{amsmath}
\usepackage{subcaption}
\usepackage{float}
\usepackage{enumitem}   
\usepackage{mathrsfs} 
\usepackage[all]{xy}
\usepackage[ruled,vlined]{algorithm2e}
\usepackage[T1]{fontenc}
\usepackage{textcomp}
\usepackage{tikz}
\usepackage{capt-of}
\usepackage{lineno}
\usepackage{mathtools}
\usepackage[colorlinks=true]{hyperref}
\hypersetup{colorlinks=true,	linkcolor=red,	filecolor=magenta, urlcolor=blue,      	citecolor=blue,}
\usepackage[capitalize,noabbrev]{cleveref}
\usepackage{etoolbox}
\usepackage{comment}
\usepackage{titlesec}
\apptocmd{\thebibliography}{\raggedright}{}{}

\newtheorem{lemma}[subsection]{Lemma}
\newtheorem{proposition}[subsection]{Proposition}
\newtheorem{theorem}[subsection]{Theorem}
\theoremstyle{definition} 
 
\newtheorem{definition}[subsection]{Definition}
\theoremstyle{remark}
\newtheorem{remark}[subsection]{Remark}
\newtheorem{corollary}[subsection]{Corollary}
\newtheorem{Problem}[subsection]{Problem}

\newtheorem{example}[subsection]{Example}

\DeclareMathOperator{\val}{val}
\DeclareMathOperator{\SL}{SL}
\DeclareMathOperator{\GL}{GL}
\DeclareMathOperator{\supp}{supp}

\DeclareMathOperator{\trop}{trop}
\DeclareMathOperator{\Sym}{Sym}
\DeclareMathOperator{\PSD}{PSD}
\DeclareMathOperator{\PD}{PD}
\DeclareMathOperator{\diag}{diag}
\DeclareMathOperator{\im}{im}

\DeclareMathOperator{\PGL}{PGL}

\newcommand{\R}{\mathbb{R}}
\newcommand{\Q}{\mathbb{Q}}
\newcommand{\Z}{\mathbb{Z}}

\newcommand{\C}{\mathbb{C}}

\newcommand{\F}{\mathbb{F}}
\newcommand{\Prob}{\mathbb{P}}

\newcommand{\K}{\mathbb{K}}

\newcommand{\Ccal}{\mathcal{C}}

\newcommand{\Scal}{\mathcal{S}}
\newcommand{\Pcal}{\mathcal{P}}
\newcommand{\Ocal}{\mathcal{O}}

\newcommand{\Lcal}{\mathcal{L}}

\newcommand\independent{\protect\mathpalette{\protect\independenT}{\perp}}
\def\independenT#1#2{\mathrel{\rlap{$#1#2$}\mkern2mu{#1#2}}}

\newcommand*\varhrulefill[1][0.4pt]{\leavevmode\leaders\hrule height#1\hfill\kern0pt}

\makeatletter
\renewcommand{\section}{\@startsection {section}{1}{\z@}%
	{-3.5ex \@plus -1ex \@minus -.2ex}%
	{2.3ex \@plus .2ex}%
	{\large \bfseries \filcenter}}
\renewcommand{\@secnumfont}{\bfseries}
\makeatother

\makeatletter
\def\blfootnote{\gdef\@thefnmark{}\@footnotetext}
\makeatother

\newcommand*{\Scale}[2][4]{\scalebox{#1}{$#2$}}%

\title{The Gaussian entropy map in valued fields}

\author{Yassine El Maazouz}

\address{Yassine El Maaouz, U.C. Berkeley, Department of statistics, 335 Evans Hall \#3860 Berkeley, CA 94720 U.S.A.}

\email{yassine.el-maazouz@berkeley.edu}

\date{October 25, 2020 }



\begin{document}

	\begin{abstract} The entropy map for multivariate real valued Gaussian distributions is the map that sends a positive definite matrix $\Sigma$ to the sequence of logarithms of its principal minors $(\log(\det(\Sigma_{I})))_I$. We exhibit the analog of this map in the non-archimedean local fields setting (like the field of $p$-adic numbers for example). As in the real case, the image of this map lies in the supermodular cone. Moreover, given a multivariate Gaussian measure on a local field, its image under the entropy map determines its pushforward under valuation. In general, this map can be defined for non-archimedian valued fields whose valuation group is an additive subgroup of the real line, and it remains supermodular. We also explicitly compute the image of this map in dimension 3.
	\end{abstract}
	
	\maketitle

	\section{Introduction and notation}\label{Section1}
	
	Gaussian measures on local fields are introduced in \cite{evans2001local}. In this text, we aim to exhibit the entropy map of these measures and discuss the properties this map satisfies. Our aim is to highlight the similarities with the real case. Before we discuss Gaussian measures on local fields (see \cref{sec:background}), we begin by reviewing the entropy map in the real setting.
	
	\subsection{Entropy of real multivariate Gaussian distributions}
	For a positive integer $d$, multivariate Gaussian distributions on $\R^d$ are determined by their mean $\mu \in \R^d$ and their positive semi-definite covariance matrix $\Sigma \in \R^{d \times d}$. Hence the natural parameter space for \emph{centered} (i.e with zero mean) Gaussian distributions on $\R^{d}$ is the positive semi-definite cone in $\R^{d\times d}$, which we denote by
	\[ 
	\PSD_d \coloneqq	 \{ \Sigma \in \Sym_d(\R), \langle x, \Sigma x \rangle \geq 0 \text{ for all } x \in \R^d  \},
	\]
	where $\Sym_d(\R)$ is the space of real symmetric matrices in $\R^{d \times d}$  and $ \langle \cdot , \cdot \rangle$ is the usual inner product on $\R^d$. \emph{Non-degenerate Gaussian} distributions are those whose covariance matrix $\Sigma$  is positive definite, i.e, $\Sigma \in \PD_d$ where
	\[
	\PD_d \coloneqq \PSD_d^{\circ} = \{ \Sigma \in  \Sym_d(\R), \langle x, \Sigma x \rangle > 0 \text{ for all non zero } x \in \R^d  \}.
	\]
	There is no shortage of instances where the PSD cone appears in probability and statistics \cite{sturmfels2010multivariate}, optimization \cite[Chapter 12]{michalek2019invitation} and combinatorics \cite{goemans1997semidefinite}.
	
	The positive definite cone has a  pleasant group-theoretic structure in the sense that its elements are in one-to-one correspondence with left cosets of  the orthogonal group $O_d(\R)$ in the general linear group $\GL_d(\R)$. The map sending the coset $AO_d(\R) \in \GL_d(\R) / O_d(\R)$  to $AA^T \in \PD_d$ is a bijection. This underscores the fact that multivariate Gaussians are tightly linked to the linearity and orthogonality structures that the Euclidean space $\R^d$ enjoys.
	
	An  important concept in statistics, probability, and information theory is the notion of entropy, which is a measure of uncertainty  and disorder  in a distribution (see \cite{entropy}). The entropy of a centered multivariate Gaussian with covariance matrix $\Sigma$ is given, up to an additive constant, by
	\[
	h(\Sigma) = - \log( | \det(\Sigma) | ) = - \log( \det(\Sigma) ) .
	\]
	
	If $X$ is a random vector in $\R^d$ with non-degenerate centered Gaussian distribution given by a covariance matrix $\Sigma \in \PD_d$, then for any subset $I$ of $[d] \coloneqq \{1,2,\dots,d\}$ the vector $X_I$ of coordinates of $X$ indexed by $I$ is also a random vector with non-degenerate Gaussian measure on $\R^{|I|}$. Moreover, its covariance matrix is $\Sigma_{I} = (\Sigma_{i,j})_{i,j \in I} \in \R^{|I| \times |I|}$, so we can define the entropy $h_I(\Sigma)$ of $X_I$ as
	\[
	h_I(\Sigma) \coloneqq h(\Sigma_I) = - \log(\det(\Sigma_{I})).
	\]
	The collection of entropy values $(h_I(\Sigma))_{I \subset [d]}$ satisfies the inequalities
	\begin{equation} \label{classic:supermod}
		h_I(\Sigma) + h_J(\Sigma) \leq h_{I \cap J}(\Sigma) + h_{I \cup J}(\Sigma) \text{ for any two subsets } I,J \subset [d].
	\end{equation}
	This is thanks to what is known as Koteljanskii's inequalities \cite{koteljanskii1963property} on the determinants of positive definite matrices, i.e,
	\begin{equation}
		\det(\Sigma_{I}) \  \det(\Sigma_{J}) \ \  \geq \ \  \det(\Sigma_{I \cap J}) \  \det(\Sigma_{I \cup J}).
	\end{equation}
	In the language of polyhedral geometry this means that the image of the entropy map
	\begin{align}\label{nonarch_entropy}
		\begin{split}
			H: \PD_d &\xrightarrow[]{} \R^{2^d}\\
			\Sigma &\mapsto  (h_I(\Sigma))_{I \subset [d]}
		\end{split}
	\end{align}
	lies inside the  \emph{supermodular} cone $\Scal_d$ in $\R^{2^d}$. This is the \emph{polyhedral cone} specified by the inequalities in (\ref{classic:supermod}), i.e,
	\[
	\Scal_d \coloneqq \{x = (x_I)_{I \subset [d]}  \in \R^{2^d}, x_{\emptyset} = 0 \text{ and } x_I + x_J \leq x_{I \cap J} + x_{I \cup J} \text{ for all } I,J \subset [d]\}.
	\]
	Since $x_{\emptyset} = 0$ for $x\in \Scal_d$ we can see $\Scal_d$ as a cone in $\R^{2^d - 1}$.
	
	\subsection{Main results}
	In this paper we deal with multivariate Gaussian distributions on local fields, and more generally non-archimedean valued fields. See Example \ref{ex:fields} for a discussion. In particular we shall define an analog to the entropy map and show that it satisfies the same set of inequalities (\ref{classic:supermod}). More precisely we prove the following:
	
	\begin{theorem}\label{thm:mainres}
		The push-forward measure of a multivariate Gaussian measure on a local field  by the valuation map is given by a tropical polynomial whose coefficients are given by the entropy  map of this measure (see \cref{thm:trop_poly}).  Moreover, these coefficients are supermodular. The entropy map is still well defined on non-archimedian valued fields in general, and remains supermodular (see \cref{thm:super_mod}.
	\end{theorem}
	
	This solves \cite[Conjecture 5.1]{maazouz2019statistics} which roughly  states that, given a multivariate Gaussian measure on a local field, its image under the entropy map determines its pushforward under valuation via a tropical polynomial. We shall break down \cref{thm:mainres} into several pieces. Namely, Theorems  \ref{thm:trop_poly}  and  \ref{thm:super_mod} for the local field case, and the discussion in Section \ref{Section4}  for the general non-archimedean valued field case.
	
	One motivation for this paper is the search for a suitable definition of  \emph{tropical Gaussian measures} \cite{tran2020tropgaussians}. Tropical stochastics has been an active research area in the recent years and has diverse applications from phylogenetics \cite{lin2018tropical,yoshida2017tropical} to game theory \cite{akian2012tropical} and economics \cite{baldwin2013tropical, tran2015product}. One appealing approach to define \emph{tropical Gaussians} is to tropicalize Gaussian measures on a valued field. Our text is organized as follows.

	In \cref{Section2} we show that tropicalizing multivariate Gaussians on local fields yields probability measures on the lattice $\Z^d$ that are determined by the entropy map via a tropical polynomial. In \cref{Section3} we show the supermodularity of the entropy map and provide a recursive algorithm to compute it. In \cref{Section4}, we explain why orthogonality is not a suitable approach to define Gaussian measures when the field $K$ is not locally compact. Nevertheless, we will see that the entropy map is still well defined and remains supermodular and we explicitly compute its image when $d = 3$.
	
	Implementations, computations and data related to this paper are made available at 
	\begin{equation}\label{algo_url}
		\text{\url{https://mathrepo.mis.mpg.de/GaussianEntropyMap/index.html}.}
	\end{equation}

	\begin{remark}
	For readers not familiar with local fields, we refer to \cite{nealKoblitz,serre2013local}. Local fields are not commonly used in statistics and probability. However, in recent years there has been a stream of literature addressing probabilistic and statistical questions in the $p$-adic setting, starting from the early work of Evans \cite{padicPhys1,evans2001local,evansNoise95} to the more recent developments \cite{ElementaryDivisors,caruso2021zeroes,kulkarni2021integGeo} to mention a few.
	\end{remark}
	
	\textbf{Acknowledgements}:  \footnote{The polyhedral geometry images in Figures \ref{fig:3d}, \ref{fig:support} and \ref{fig:Ccal_Pcal} were drawn using Polymake \cite{polymake:2000}.}
	 The author would like to thank the Max Planck Institute for Mathematics in the Sciences for the generous hospitality while working on this project. He would also like to thank Bernd Sturmfels and Ian Le for  valuable mathematical discussions. The author is grateful to Avinash Kulkarni for the numerous and valuable exchanges while writing this paper. Many thanks also to Rida Ait El Mansour and Adam Quinn Jaffe  for their remarks on early drafts of this manuscript. Finally, the author thanks the anonymous referee for valuable comments and remarks.

	\section{Background on valued fields and Gaussian measures}\label{sec:background}
		
		This section is meant to collect the basic facts and result that we will need in our discussion. Most of these results can be found in the literature on valued fields in number theory \cite{serre2013local,weil2013basic,engler2005valued} and functional analysis \cite{van1978non,schikhof1984ultrametric,schikhof2007ultrametric}.
		
		\subsection{Valued fields}
		Let $K$ be a field with an \emph{additive  non-archimedean valuation } $\val: K \xrightarrow[]{} \R \cup \{+ \infty\}$  with valuation group $\Gamma \coloneqq \val(K^\times)$. The valuation map $\val$ defines an equivalence class of \emph{exponential valuations} or \emph{absolute values}  $|\cdot|$ on $K$ via $|x| \coloneqq a^{-\val(x)}$ (where $a \in (1,\infty)$) and hence also a topology on $K$. The valuation $\val$ is called \emph{discrete} if its valuation group $\Gamma$ is a discrete subgroup of $\R$ which, by scaling $\val$ suitably, we can always assume to be $\Z$ (we then call $\val$ a \emph{normalized valuation}). In the discrete valuation case we fix a uniformizer $\pi$ of $K$, i.e, an element $\pi \in K$ with $\val(\pi) = 1$. We  denote by $\Ocal \coloneqq \{ x \in K, \val(x) \geq 0\}$ the valuation ring of $K$; this is a local ring with unique maximal ideal $\mathfrak{m} \coloneqq \{ x \in K, \val(x) > 0\}$ and residue field $k \coloneqq \Ocal / \mathfrak{m}$. When the valuation is discrete, the ideal $\mathfrak{m}$ is generated in $\Ocal$ by $\pi$ i.e $\mathfrak{m} \coloneqq \pi \Ocal$. We mention typical examples of such fields in \cref{ex:fields}.
		
		\begin{example} \label{ex:fields}
			\begin{enumerate}[label=(\arabic*)]
				\setlength\itemsep{0.1em} 
				
				\item The field $\F_q((t))$ of Laurent series in one variable  with coefficients in the finite field $\F_q$.
				\item The fields $\R((t))$ or $\C((t))$ of Laurent series with complex or real coefficients. These are fields with an infinite residue field but still in discrete valuation $\Gamma = \Z$.
				\item The fields $\R\{\{t\}\} = \cup_{n \ge 1} \R((t^{1/n}))$ and $\C\{\{t\}\} = \cup_{n \ge 1} \C((t^{1/n}))$  of Puiseux series in $t$. In this case the valuation group $\Gamma = \Q$ is dense.
				\item Another interesting field is the field of generalized Puiseux series $\K$ which has valuation group $\Gamma = \R$. This field consists of formal series $\bm{f} = \sum_{\alpha \in \R} a_\alpha t^\alpha$ where $\supp(\bm{f}) \coloneqq \{\alpha \in \R : a_\alpha \neq 0\} $ is either finite or has $+ \infty$ as the only accumulation point. See \cite{allamigeon2018log} and references therein.
				\item All the previous fields have the same characteristic as their residue fields.  Interesting examples in mixed characteristic are the field of $p$-adic numbers $\Q_p$ where $p$ is prime, its algebraic closure $\overline{\Q}_p$ and the field of $p$-adic complex numbers $\C_p$ (completion of $\overline{\Q}_p$).
			\end{enumerate}
		\end{example}
		
\subsection{Local fields} These are valued fields that are locally compact. In this section let us assume that $K$ is locally compact. It is then known that  $K$ is isomorphic to a finite field extension of $\Q_p$ or $\F_q((t))$ and that its valuation group $\Gamma$ is discrete in $\R$, and its residue field $k$ is finite. In this case, by convention, the absolute valued on $K$ is defined as $|x| = q^{-\val(x)}$ (so we choose $a = q$), and there exist a unique Haar measure $\mu$ on $K$ such that $\mu(\Ocal) = 1$.

\subsection{Lattices} 	Let $d \geq 1$ an integer. We call a \emph{lattice} in $K^d$ any $\Ocal$-submodule $\Lambda \coloneqq \bigoplus_{i = 1}^{n} \Ocal a_i$ generated by a basis $(a_1, \dots, a_d)$ of $K^d$. The basis $(a_1,\dots, a_d)$ that generates $\Lambda$ is not unique. We can write $\Lambda = A \Ocal^d$ where $A$ is the matrix with columns $a_1, \dots, a_d$, which is then called a \emph{representative} of $\Lambda$. The elements $U$ of the group $\GL_d(K)$ that leave $\Ocal^d$ invariant (i.e $U \Ocal^d = \Ocal^d$) are exactly the matrices $U \in \GL_d(\Ocal)$ with entries in $\Ocal$ whose inverse has all entries in $\Ocal$. The group $\GL_d(\Ocal)$ then plays the role of the orthogonal group $O_d(\R)$ \cite[Theorem 2.4]{evans2019rotatable}. Then, like positive definite matrices matrices, lattices are  in a one-to-one correspondence with left cosets $\GL_d(K)/ \GL_d(\Ocal)$, in particular, any two representatives of a lattice $\Lambda$ are elements of the same left coset. A lattice $\Lambda$ is called \emph{diagonal}\footnote{Homotethy classes of diagonal lattices form what is called an \emph{apartment} in the theory of buildings.} if it admits a diagonal matrix as a representative. Let us now state a result on lattices over valued fields that will be useful in our discussion.

\begin{lemma}\label{lem:G-inv}
	For any  two lattices $\Lambda, \Lambda'$ there exists an element $g \in \GL_d(K)$ such that $g. \Lambda$ and $g .\Lambda '$ are both diagonal lattices. \footnote{This is in fact a property of buildings: any two chambers belong to a common apartment. See \cite{abramenko2008buildings}.}
\end{lemma}
\begin{proof}
	It suffices to show this when $\Lambda$ is the standard lattice $\Lambda = \Ocal^d$. Let $A \in \GL_d(K)$ be a representative of $\Lambda'$. Thanks to the non-archimedean single value decomposition (see \cite[Theorem 3.1]{evans2002elementary}), there exists a diagonal matrix $D \in \GL_d(K)$ and  $U,V \in \GL_d(\Ocal)$ such that $A = U D V$. Hence we deduce that $\Lambda' = UD \Ocal^d$. Picking $g = U^{-1}$ yields $g \Lambda = U^{-1} \Ocal^d = \Ocal^d  $ and $g \Lambda' = D\Ocal^d$.
\end{proof}

\subsection{Gaussian measures} Suppose that $K$ is a local field and $d$ is a positive integer. As shown by Evans \cite{evans2001local},  one can define multivariate  \emph{Gaussian} measures on $K^d$ using non-archimedean orthogonality. It turns out that these measures are precisely the uniform distributions on $\Ocal$-submodules of $K^d$. The non-degenerate Gaussians on $K^d$ are then parameterized by full rank submodules of $K^d$ i.e. lattices.

 For a lattice $\Lambda$ in $K^d$ we denote by $\Prob_\Lambda$ the Gaussian measure on $K^d$ given by  $\Lambda$, i.e. the uniform probability measure on $\Lambda$. If $f_\Lambda$ denote the density (with respect to the Haar measure $\mu^{\otimes d}$) of $\Prob_\Lambda$, then 
  \[
	  f_\Lambda(x) = \bm{1}_\Lambda(x) / \mu^{\otimes d }(\Lambda), \quad x \in K^d,
  \]
  where $\bm{1}_\Lambda$ is the set indicator function of $\Lambda$.

One can then think of lattices as an analogues for the positive definite covariance matrices in the real case since they parametrize non-degenerate multivariate Gaussian measures. In the language of group theorists, one can think of the Bruhat-Tits building for the reductive group $\PGL_d(K)$ \cite{abramenko2008buildings} as the parameter space for non-degenerate Gaussians up to scalar multiplication.

	\section{The entropy map of local field Gaussian distributions }\label{Section2}

	In this section we assume that $K$ is a local field and we fix a positive integer $d \geq 1$ and a lattice $\Lambda$ in $K^d$. We recall that there is a unique Haar measure $\mu^{\otimes d}$ on $K^d$  which is the product measure induced by $\mu$ on $K$. Letting $A$ be a representative of the lattice $\Lambda$, i.e. $\Lambda = A \Ocal^d$, we can define the \emph{entropy} $h(\Lambda)$ of the lattice $\Lambda$ as 
	\[
	h(\Lambda) = \val(\det(A)).
	\] 
	This is a well defined quantity since any other representative of $\Lambda$ is of the form $A U$ where $U \in \GL_d(\Ocal)$ and $\det(U) \in \Ocal ^\times$ is a unit, so $\val(\det(U)) = 0$. This definition lines up with the definition in the real case because $\val(x) = -\log_q(|x|)$ where $|\cdot|$ is the absolute value on $K$, so we get 
	\[
	h(\Lambda) = \val(\det(A)) = -\log_q(| \det(A) | ) .
	\]
	
	The following proposition justifies the nomenclature ``entropy'' and relates the entropy $h(\Lambda)$ of a lattice $\Lambda$ to its measure $\mu^{\otimes d}(\Lambda)$.

	\begin{proposition}\label{prop:entropyandMeasure}
		We have $\mu^{\otimes d}(\Lambda) = q^{ - h(\Lambda)}$. Moreover, the quantity $h(\Lambda)$ is the differential entropy of the Gaussian measure $\Prob_\Lambda$, i.e, 
		\[
		h(\Lambda) = \int_{K^d} \log_q( f_\Lambda(x)) \Prob_\Lambda (dx) .
		\]
	\end{proposition}
	\begin{proof}
			Let $A$ be a representative of $\Lambda$. Thanks to the non-archimedean single value decomposition (see \cite[Theorem 3.1]{evans2002elementary}), we can write $A = U D V$, where $U,V \in \GL_d(\Ocal)$ are two orthogonal matrices and $D$ is a diagonal matrix.	Then we have $\Lambda = U D . \Ocal^d $. Since orthogonal linear transformation in $K^d$ preserve the measure, we have $\mu^{\otimes d}(\Lambda) = \mu^{\otimes}(D.\Ocal^d)$. Let $\alpha_1,\dots, \alpha_d$ be the diagonal entries of $D$. Then we have $\mu^{\otimes d}( \Lambda) = \mu^{\otimes d} ( \bigoplus_{i = 1}^{d} \alpha_i \Ocal  ) = q^{- \val(\alpha_1) - \dots - \val(\alpha_d)}$. But $\val(\alpha_1) + \dots + \val(\alpha_d) = \val(\det(A)) = h(\Lambda)$. The second statement follows from the immediate computation:
		\begin{align*}
			\int_{K^d} \log_q( f_\Lambda(x)) \Prob_\Lambda (dx)  = \int_{K^d} \log_q( f_\Lambda(x))   f_\Lambda(x) \mu^{\otimes d}(dx) =  h(\Lambda).
		\end{align*}
	\end{proof}
	
	For a subset $I$ of $[d] \coloneqq \{1,2,\dots, d\}$ we denote by $\Lambda_I$ the image of $\Lambda$ under the projection onto the space $K^{|I|}$ of coordinates indexed by $I$. This is also a lattice in the space $K^{|I|}$.  So, for any subset $I \subset [d]$, we can define the entropy $h_I(\Lambda)$ of the lattice $\Lambda_I$. We can then define the entropy map 
	\begin{align*}
		H: \GL_d(K)/\GL(\Ocal) &\xrightarrow[]{} \R^{2^d}\\
		\Lambda &\mapsto  (h_I(\Lambda))_{I \subset [d]}
	\end{align*}
	where $h_{\emptyset}(\Sigma) = 0$ by convention.	If $A$ is a representative of $\Lambda$ with columns $a_{1}, \dots, a_{d}$, then the lattice $ \Lambda_I$ is the lattice generated over $\Ocal$ by the vectors $a_{i,I}$ which are the sub-vectors of the $a_i$'s with coordinates indexed by $I$. So we can compute $h_I(\Lambda)$ from the matrix $A$ by
	\begin{equation} \label{entropycoefs}
		h_I(\Lambda) = \min_{J \subset [d], |J| = |I|}   \val(\det(A_{I \times J})) ,
	\end{equation}
	where $A_{I \times J}$ is the matrix extracted from $A$ by taking the rows indexed by $I$ and the columns indexed by $J$,  i.e.  $A_{I \times J} = (A_{i,j})_{i \in I, j\in J}$.
	
	Now let $X$ be a $K^d$-valued random variable with Gaussian distribution $\Prob_\Lambda$ given by $\Lambda$. So for any measurable set $B$ in the Borel $\sigma$-algebra of $K^d$,
	\[
	\Prob_{\Lambda}(X \in B) = \frac{\mu^{\otimes d} ( \Lambda \cap B)}{\mu^{\otimes d} ( \Lambda)},
	\]
	and $V \coloneqq \val(X)$ its image under coordinate-wise valuation. Notice that, since $\Prob_\Lambda(X_i = 0) = 0$ for any $i \in \{1,\dots,d\}$, the vector $V$ is almost surely in $\Z^d$ . By definition the distribution of $V$ is the push-forward of the distribution of $X$ by the map $\val$. We are interested in the distribution of the valuation vector $V$ and to determine it we compute its \emph{tail distribution function} $Q_\Lambda$ which is defined on $\R^d$ as
	\[
	Q_\Lambda(v) \coloneqq   \Prob_{\Lambda}( V \geq v) \text{ for any  } v \in \R^d,
	\] 
	where $\geq$ is the coordinate-wise partial order on $\R^d$. Since $V$ takes values in $\Z^d$ this, function is completely determined by its values for $v \in \Z^d$. For a vector $v=(v_1, \dots, v_d) \in \Z^d$ let us denote by $\bm{\pi^v}$ the $\Ocal$-module generated by the basis $\pi^{v_i} e_i$ where $e_1,\dots,e_d$ is the standard basis of $K^d$ i.e.
	\[
		\bm{\pi^v} = \pi^{v_1} \Ocal e_1 \oplus \dots \oplus \pi^{v_d} \Ocal e_d.
	\]

		\begin{definition}\label{def:phi}
		We define the \emph{logarithmic tail distribution function}  $\varphi_{\Lambda}$ as
		\[
		\varphi_\Lambda \colon \Z^d \to \Z , \quad v \mapsto -\log_q(Q_\Lambda(v)).
		\]
		\end{definition}
		The following lemma relates the tail distribution function $\varphi_\Lambda$ with the entropy $h(\Lambda)$ of the lattice $\Lambda$.
		
		\begin{lemma} \label{lem:phiandEntropyandGroupindex}
			We have $\varphi_{\Lambda}(v) = h(\Lambda \cap \bm{\pi^v}) - h(\Lambda)$. Moreover, if  $[\Lambda : \Lambda \cap \bm{\pi^v}]$ denotes the index of $\Lambda \cap \bm{\pi^v}$ as a subgroup of $\Lambda$ then we also have 
			\[
				Q_{\Lambda}(v) = 1 / [ \Lambda : \Lambda \cap \bm{\pi^v}].
			\]
		
		\end{lemma}
		\begin{proof}
			By definition we have $Q_\Lambda(v) = \Prob_\Lambda(X \in \bm{\pi^{v}} ) = \frac{\mu^{\otimes d}(\bm{\pi^{v}} \cap \Lambda )}{\mu^{\otimes d}(\Lambda)} $. So by virtue of \cref{prop:entropyandMeasure} we deduce that $Q_\Lambda(v) = q^{h(\Lambda) - h(\Lambda \cap \bm{\pi^v}) }$. The first statement then follows from the definition of $\varphi_\Lambda$ (\cref{def:phi}). For the second statement, by definition, $\Lambda$ can be partitioned into $[ \Lambda : \Lambda \cap \bm{\pi^v}]$ cosets of $\Lambda \cap \bm{\pi^v}$. Since the Haar measure $\mu^{\otimes d}$ is translation invariant all of these cosets have the same measure i.e. $\mu^{\otimes d}(\Lambda) =  [ \Lambda : \Lambda \cap \bm{\pi^v}] \mu^{\otimes d}(\Lambda \cap \bm{\pi^v})$. The result then follows from the fact that $Q_\Lambda(v) =  \frac{\mu^{\otimes d}(\bm{\pi^{v}} \cap \Lambda )}{\mu^{\otimes d}(\Lambda)} $ and \cref{def:phi}.
		\end{proof}

	 Next, we introduce a technical tool that we will be using in the proof of our first result.
	 
	\begin{definition}
		For any $\ell \in \{0,\dots,d\}$ we define the $\ell$\emph{-distance}  $\phi_\ell(\Lambda, \Lambda')$  of two lattices $\Lambda, \Lambda'$ as the minimum of $\val(\det(x_1,\dots,x_\ell,y_1,\dots,y_k))$ among all possible choices of $x_1,\dots,x_{\ell} \in \Lambda $ and $y_1,\dots, y_k \in \Lambda' $ where $k = d-\ell$. 
	\end{definition}
	Since for any $g \in \GL_d(K)$ ,  $x_1, \dots x_{\ell} \in \Lambda$  and $y_1, \dots, y_k \in \Lambda'$ we have
	\[
	\Scale[0.95]{\val(\det(g x_1, \dots g x_{\ell}, g y_1, \dots , g y_k)) = \val(\det(x_1,\dots,x_\ell,y_1,\dots,y_k)) + \val(\det(g)),}
	\]
	we can see that $\phi_\ell$ satisfies the following property:
	\[
	\phi_\ell(g.\Lambda, g.\Lambda') = \phi_\ell(\Lambda, \Lambda') + \val(\det(g)).
	\]
	We then deduce that the quantity $\phi_\ell(\Lambda,\Lambda') - h(\Lambda' ) $ is invariant under the action $\GL_d(K)$, i.e, for any $g \in \GL_d(K)$ we have
	\[
	\phi_\ell(g.\Lambda,g.\Lambda') - h(g.\Lambda' ) = \phi_\ell(\Lambda,\Lambda') - h(\Lambda' ).
	\]
	
	When the second lattice $\Lambda' = \bm{\pi^v}$ is diagonal and $\Lambda$ has representative $A \in \GL_d(K)$, the optimal choice for the vectors $x_1,\dots ,x_\ell$ and $y_1,\dots,y_k$ is when the vectors $x_1,\dots,x_\ell$ are among the columns $a_1,\dots,a_d$ of $A$   and the vectors $y_1,\dots,y_k$ are among the vectors $\pi^{v_i} e_i$ where $(e_i)_{1\leq i \leq d}$ is the standard basis of $K^d$. So we deduce that $\phi_\ell(\Lambda, \bm{\pi^v})$ can be computed as follows:
	\[
	\phi_\ell(\Lambda, \bm{\pi^v}) = \min_{  \substack{ I,J \subset [d]\\  |I| = |J| = \ell}} \left( \val(\det(A_{I \times J}))+ \sum_{ j \not \in J} v_j \right).
	\]
	So we also get 
	\begin{equation}\label{entropy_polyn}
		\phi_\ell(\Lambda, \bm{\pi^v}) - h(\bm{\pi^v}) = \min_{ \substack{ I,J \subset [d]\\  |I| = |J| = \ell}}\left( \val(\det(A_{I \times J}) ) - \sum_{ j  \in J} v_j \right).
	\end{equation}
	In the special case $\Lambda = \bm{\pi^a} $, for $a \in \Z^d$, the determinant of $A_{I\times J}$ in the above optimization problem is $0$ whenever $J \neq I$, since we can choose $A$ to be diagonal. So we get the following 
	\[
	\phi_\ell(\bm{\pi^{a}}, \bm{\pi^v}) - h(\bm{\pi^v}) = \min_{I  \subset [d], |I| = \ell} \left( \sum_{i \in I} a_i - \sum_{ i  \in I} v_i \right).
	\]
	
	\begin{theorem}\label{thm:trop_poly}
		The logarithmic tail distribution function $\varphi_{\Lambda}$ is a tropical polynomial on $\Z^d$ given by
		\begin{equation}\label{trop_poly}
			\varphi_\Lambda(v) = \max_{I \subset [d]} (v_I - h_I(\Lambda)).
		\end{equation}
		
	\end{theorem}

	\begin{proof}
		First we show this for a diagonal lattice $\Lambda = \bm{\pi^{a}} $ where $a \in \Z^d$. For any $v \in \Z^d$, let $a\vee v$ the vector with coordinates $\max(a_i,v_i)$.  We have $\bm{\pi^a} \cap \bm{\pi^v} = \bm{\pi^{a \vee v}}$ so we get the entropy $h(\bm{\pi^a}) = \sum_{i=1}^{d} a_i $ and $h( \bm{\pi^a} \cap \bm{\pi^v} ) = h(  \bm{\pi^{a \vee v}} ) = \sum_{i=1}^{d} \max(a_i,v_i)$. Hence we have 
		\[
		\varphi_{\Lambda}(v) =  h( \bm{\pi^{a}} \cap \bm{\pi^v} )  - h( \bm{\pi^{a}} ) = \max_{I \subset [d]} \left( \sum_{i \in I} v_i + \sum_{i \not \in I} a_i  \right) - \sum_{i = 1}^{d} a_i  = \max_{I \subset [d]} (v_I - a_I ),
		\]
		and $h_I(\bm{\pi^a}) = a_I$. So the theorem holds for diagonal lattices. To see why it also holds for a general lattice $\Lambda$, first notice that in the diagonal case $\Lambda = \bm{\pi^a}$ we have
		\[		
		\varphi_{\Lambda}(v) = - \min_{ \ell = 0,\dots, d} \left( \phi_\ell( \Lambda, \bm{\pi^v}) - h(\bm{\pi^v}) \right).
		\]
		Secondly, notice that the right hand side of the previous equation is invariant under the action of $\GL_d(K)$. So for $g \in \GL_d(K)$, 
		\[
		\min_{ \ell = 0,\dots, d} \left( \phi_\ell( g .\Lambda, g . \bm{\pi^v}) - h(g . \bm{\pi^v}) \right) =  \min_{ \ell = 0,\dots, d} \left( \phi_\ell( \Lambda, \bm{\pi^v}) - h(\bm{\pi^v}) \right).
		\]
		By  \cref{lem:phiandEntropyandGroupindex}, we have $\varphi_\Lambda(v) = \log_q([\Lambda : \Lambda \cap \bm{\pi^v}] ) = \log_q ([g .\Lambda : g. \Lambda \cap g.\bm{\pi^v}])$. Now fix a general lattice $\Lambda$ and $v \in \Z^d$. Also, by Lemma \ref{lem:G-inv}, there exists $g \in \GL_d(K)$ such that $g \Lambda$ and $g \bm{\pi^v}$ are both diagonal, so
		\begin{align*}
			\varphi_\Lambda(v) = \log_q( [g .\Lambda : g. \Lambda \cap g.\bm{\pi^v}] ) &=  - \min_{ \ell = 0,\dots, d} \left( \phi_\ell( g .\Lambda, g . \bm{\pi^v}) - h(g . \bm{\pi^v}) \right)\\
			& =  - \min_{ \ell = 0,\dots, d} \left( \phi_\ell( \Lambda,   \bm{\pi^v}) - h( \bm{\pi^v}) \right) .
		\end{align*}
		Hence, we deduce, thanks to equation (\ref{entropy_polyn}), that
		\[
		\varphi_{\Lambda}(v) = - \min_{\ell = 0, \dots, d} \left(  \min_{  \substack{ I,J \subset [d]\\  |I| = |J| = \ell} } \left( \val(\det(A_{I \times J} ) )  - \sum_{ j  \in J} v_j \right)  \right) .
		\]
		We can simplify this thanks to equation (\ref{entropycoefs}) to get the desired equation (\ref{trop_poly}).
	\end{proof}

	So the distribution of the random vector of valuations $V$ is given by a tropical polynomial $\varphi_{\Lambda}$ via its tail distribution function $Q_\Lambda$. The coefficients of this polynomial are exactly the entropies $h_I(\Lambda)$. Now we prove a couple of interesting properties of $\varphi_\Lambda$, namely how the coefficients $h_I(\Lambda)$ behave under diagonal scaling and permutation of coordinates of the random vector $X$. To this end, let us denote by $D_a = \diag(a_1,\dots, a_n)$ the diagonal matrix with coefficients $a_i \in K$ and $P^\sigma$ the permutation matrix corresponding to a permutation $\sigma$ of $[d]$ i.e $ P^\sigma_{i,j} = 1 $ when $j = \sigma(i)$ and $0$ otherwise.

	\begin{lemma} \label{lem: diag_perm}
		Let $\Lambda$ be a lattice in $K^d$,  $a \in K^d$ and $\sigma$ a permutation of $[d]$. We have the following:  
		
		\begin{center}
			$ h_I(D_a \Lambda)  = h_I(\Lambda) + \sum\limits_{i \in I} \val(a_i) \text{ and } h_I( P^\sigma \Lambda) =  h_{\sigma(I)}(\Lambda)$.
		\end{center}
	\end{lemma}
	
	\begin{proof}
		For $I \subset [d]$, we have $ h_I (D_a \Lambda) = \min\limits_{ |J| = |I|  } \val(\det((D_aA)_{I\times J})) $, where $A$ is any representative of $\Lambda$. Since all the lines of $D_aA$ are multiples of those of $A$ by the scalars $a_i$ we deduce that $\det((D_a A)_{I\times J}) = \det(A_{I\times J}) \prod_{i \in I} a_i$ and hence we get 
		\[
		h_I(D_a \Lambda)  = h_I(\Lambda) + \sum_{i \in I} \val(a_i).
		\]
		Similarly we can see the effect the permutation of coordinates of $X$ has on the vector of entropies $ H(\Lambda) = ( h_I(\Lambda) )_{I \subset [d]}$.
	\end{proof}

	\section{Supermodularity of the entropy map}\label{Section3}

	As it is the case for real Gaussians, we would like the vector of entropies $ H(\Lambda) : = (h_I(\Lambda))$ to have values in the supermodular cone $\Scal_d$ as conjectured in \cite{maazouz2019statistics}. As a first step towards proving this result, notice that the previous lemma implies that if $\Lambda$ is a lattice such that $H(\Lambda) \in \Scal_d$, then for any diagonal matrix $D_a$ we still have $H(D_a \Lambda ) \in \Scal_d $ and $H( P^{\sigma} \Lambda) \in \Scal_d$ for any permutation $\sigma$ of $\{1 , \dots , d\}$.
	
	\begin{definition}[Hermite normal form  \footnote{The curious reader can see \cite[Chapter II]{weil2013basic} and \cite[Proposition 4.2]{maazouz2019statistics} for more details.}]\label{Hermite}
		Every lattice $\Lambda$ in $K^d$ has a representative $A$ in \emph{Hermite normal form}, i.e, a matrix $A= (A_{ij})$ in $\GL_d(K)$ satisfying the following conditions:
		\begin{enumerate}[label=(\roman*)] 
			\item $A$ is lower triangular i.e. $A_{ij} = 0$ whenever $i < j$. 
			\item For any $1\leq j < i\leq d$ we have either $\val(A_{ij}) < \val(A_{jj})$ or $A_{ij} = 0$.
			\item The diagonal coefficients $A_{ii}$ are of the form $A_{ii} = \pi^{a_i}$ for some $a_i \in \Z$.
		\end{enumerate}	
	\end{definition}		 
	
	Now we can state the second result of this section concerning the supermodularity of the entropy map. But, before we do that, we give an equivalent definition of the supermodular cone as follows:
	\[
	\Scale[0.95]{\Scal_d = \left\{ (x_I)_{I \subset [d]} \in \R^{2^d}\colon  \begin{cases} x_{\emptyset } = 0 \\  x_{Ii} + x_{Ij} \leq x_{I} + x_{Iij},   \text{ for any } I \subset [d], i\neq j \in [d] \setminus I  \end{cases}     \right\}}
	\]
	where we write $Ii$ instead of $I \cup \{i \}$. These are the facet-defining inequalities of the cone $\Scal_d$ and there are $d(d-1) 2^{d-3}$ of them. See \cite{kuipers2010generalization} and references therein.
	
	\begin{theorem}	 \label{thm:super_mod}
		The image of the map $H : \Lambda \to ( h_I(\Lambda)) _{I \subset[d]}$ lies in the supermodular cone $\Scal_d$, i.e,  for any subset $I \subset [d]$ with $|I| \leq d-2 $ and $i \neq j \in [d] \setminus I$,
		
		\[
		h_{I i }(\Lambda) + h_{I j}(\Lambda) \leq h_{I}(\Lambda) + h_{Iij}(\Lambda).
		\]
	\end{theorem}
	\begin{proof}[ Proof]
		We prove this by induction on $d$. The result is trivial for $d =1,2$. Assume that it holds for lattices in $K^r$ for any $ r \leq d$, where $d \geq 3$. Let $\Lambda$ be a lattice in $K^{d}$ and $A$ its Hermite normal form. For any $I \subset [d]$ of size $|I| < d-2$ the inequality  $h_{I i }(\Lambda) + h_{I j}(\Lambda) \leq h_{I}(\Lambda) + h_{Iij}(\Lambda)$ holds for any $i \neq j$ not in $I$ thanks to the induction hypothesis. This is because, when $|I| \leq d-2$, we are working on the lattice $\Lambda_{Iij}$ which is a lattice in dimension less than $d$. Then, it suffices to show the inequality when $I$ has size $d-2$. By Lemma \ref{lem: diag_perm} we can assume that $I = \{1,\dots, d-2\}$ and $i = d-1$ and $j = d$ (if not, we can just act on $\Lambda$ by a suitable permutation matrix). Let us write down the matrix $A$ as follows
		
		\[
		A = \begin{pmatrix}
			\pi^{a_1}   &            0           &    \dots             &               0           &          0               &      0          \\
			\ast          &       \pi^{a_2}     &    \ddots           &             \vdots      &     \vdots            &     \vdots    \\   
			\vdots      &       \ddots         &    \ddots           &               0           &              0           &    0  	      \\      
			\ast          &       \dots          &          \ast         &       \pi^{a_{d-2}}   &     0                    &   0             \\     
			\ast          &       \dots          &          \ast         &               \ast       &     \pi^{a_{d-1}}    &   0              \\       
			\ast          &       \dots          &          \ast         &   \ast                   &  x                      &  \pi^{a_{d}}   \\       
		\end{pmatrix}.
		\]
		Recall that since $A$ is the Hermite form of $\Lambda$ we have $\val(x) < a_d $ or $x = 0$. Now we have
		\begin{align*}
		h_{Ii}(\Lambda) = a_1 + \dots + a_{d-1}  &  , \quad    \quad             &   h_{Ij}(\Lambda) = a_1 + \dots +a_{d-2} + \min(\val(x), a_d) \\
		h_{I}(\Lambda) = a_1 + \dots  + a_{d-2}    &  , \quad    \text{ and }     &   h_{Iij}(\Lambda) = a_1 + \dots +a_{d}.
		\end{align*}
		The inequality  $h_{I i }(\Lambda) + h_{I j}(\Lambda) \leq h_{I}(\Lambda) + h_{Iij}(\Lambda)$ then holds simply because $\min(\val(x), a_d) \leq a_d$ and this finishes the proof.		
	\end{proof}
	
	This theorem underlines another similarity between the local field Gaussians defined in \cite{evans2001local} and classical multivariate Gaussian measures. From Lemma (\ref{lem: diag_perm}) we can see that acting on $\Lambda$ by a diagonal matrix just moves the point $H(\Lambda) \in \Scal_d$ in parallel to the \emph{lineality} space of the cone $\Scal_d$, that is, the biggest vector space contained in $\Scal_d$.

	The classical entropy map is tightly related to conditional independence. More precisely, if $\Sigma \in \PD_d$ and $X$ is a Gaussian vector with covariance matrix $\Sigma$, then for any $I \subset [d]$ and $i \neq j$ not in $I$ the variables $X_i$ and $X_j$ are independent given the vector $X_I$ if and only if $h_{Ii} ( \Sigma ) + h_{Ij}(\Sigma) = h_I(\Sigma) + h_{Iij}(\Sigma)$ and we write
	\[
	X_i \independent X_j | X_I \iff h_{Ii} ( \Sigma ) + h_{Ij}(\Sigma) = h_I(\Sigma) + h_{Iij}(\Sigma).
	\]
	This means that the conditional independence models are exactly the inverse images by $H$ of the faces of $\Scal_d$ \cite[Proposition 4.1]{sturmfels2009open}.	It turns out that, in the local field setting, the non-archimedian entropy map $H$ defined in (\ref{nonarch_entropy}) also encodes conditional independence information on the coordinates of the random Gaussian vector $X$ as stated in the following proposition.
	
	\begin{proposition}
		Assume $d \ge 2$ and let $I$ be a subset of $[d]$ and $i\neq j \in [d] \setminus I$ two distinct integers. Let $\Lambda$ be a lattice in $K^d$ and $X$ a random Gaussian vector with distribution given by $\Lambda$. Then the conditional independence statement $X_i \independent X_j | X_I$ holds if and only if $h_{I i}(\Lambda) + h_{I j}(\Lambda) = h_{I}(\Lambda) + h_{I i j}(\Lambda)$.
	\end{proposition}
	\begin{proof}
		Using Lemma \ref{lem: diag_perm} we reduce to the case $I = [r]$ where $r \leq d-2$ ,  $i = r+1$ and $j = i+1$. Let $A = (a_{i,j})$ be the unique representative in Hermite form of $\Lambda$. We claim that  $X_{i} \independent X_{j} | X_{I}$ if and only if $a_{j,i} = 0$. To see why, let $Z = A^{-1} X$ which is a Gaussian vector whose distribution is the uniform on $\Ocal^d$. We have $X_{i} = a_{i,1} Z_1 + \dots + a_{i,i} Z_{i} $ and $X_{j} = a_{j,1} Z_1 + \dots + a_{j,j} Z_{j} $. Since $Z_I = A_{I,I}^{-1} X_I$, given $X_I$ we know $Z_I$ and vice-versa. Hence $X_{i} \independent X_{j} | X_{I}$ holds if and only if $ (a_{j, i} Z_{i} + a_{j,j} Z_{j}) \independent Z_{i}$. This happens if and only if the vectors $(1,0)$ and $(a_{j,i} , a_{j, j})$ in $K^2$ are orthogonal (see \cite{evans2001local}). This is equivalent to $\val(a_{j,j}) \leq \val(a_{j,i})$ which means that $a_{j,i} = 0$ since $A$ is in Hermite form. On the other hand, since $A$ is lower triangular, we have the following
		\begin{align*}
	\Scale[0.94]{h_{I}(\Lambda) = \val(\det(A_{I\times I}))} \ &, \ \Scale[0.94]{ h_{I i}(\Lambda) = h_{I}(\Lambda) + \val(a_{i,i}) } \\
	\Scale[0.94]{  h_{I j}(\Lambda) = h_{I}(\Lambda) + \min( \val(a_{j,i}) , \val(a_{j,j}) )  }  & \Scale[0.94]{ \text{ and  }  h_{I i j}(\Lambda) = h_{I}(\Lambda) + \val(a_{i,i}) + \val(a_{j,j}). }
		\end{align*}
		So the equality $ h_{I i}(\Lambda) + h_{I j}(\Lambda) = h_{I}(\Lambda) + h_{I ij}(\Lambda)  $ holds if and only if $ \val(a_{j,j}) \leq \val(a_{j,i}) $ since $A$ is the Hermite form of $\Lambda$ this happens if and only if $a_{j,i} = 0$. In combination with the calculation above,  this finishes the proof.
	\end{proof}
	In other terms, the conditional independence statement  $X_i \independent X_j | X_I$ holds if and only if the entropy vector $H(\Lambda) = (h_I(\Lambda))$ is on the face of the polyhedral cone $\Scal_d$ cut by the equation $h_{I i} (\Lambda) + h_{I j} (\Lambda) = h_{I } (\Lambda) + h_{I i j} (\Lambda) $. This gives an analogue of \cite[Proposition 4.1]{sturmfels2009open}.
	\begin{corollary}
		The Gaussian conditional independence models are exactly those subsets of lattices that arise as inverse images of the faces of $\Scal_d$ under the map $H$.
	\end{corollary}
	\begin{proof}
		Follows immediately from the previous proposition.
	\end{proof}
	This underlines the importance of the map $H$, and also gives reason to think that the suitable analogue of the positive definite cone on local fields is the set of lattices or more precisely the Bruhat-Tits building \cite{abramenko2008buildings,maazouz2019statistics}. A  hard question in information theory  for classical multivariate Gaussians is to describe the image of the entropy map \cite{sturmfels2009open}. This problem turns out to be difficult in this setting as well.
	\begin{Problem}\label{problem}
		Characterize the image of the entropy map $H$ and describe how it intersects the faces of $\Scal_d$. What can you say about the fibers of this map?
	\end{Problem}
	\begin{remark}\label{rem:lineality}
		We recall that for any $d \geq 1$ the image  $\im(H)$ is invariant under the action of the symmetric group and by translation in parallel to the lineality space of $\Scal_d$. This is thanks to Lemma \ref{lem: diag_perm}. We will provide an answer for Problem \ref{problem} when $d = 2, 3$ in the end of Section \ref{Section4}.
	\end{remark}
	We now provide an algorithm to compute the entropy vector $H(\Lambda)$, i.e, the coefficients of the polynomial $\varphi_{\Lambda}$. This relies on computing the Hermite form rather than directly solving the optimization problems given by equation (\ref{entropycoefs}).
	
	\begin{algorithm}[H] \label{algo}
		\SetAlgoLined
		\KwIn{ A full rank matrix  $A = (a_1, \dots , a_n) \in K^{d \times n}$ with $n \geq d$  generating $\Lambda$}
		\KwOut{The entropy vector $H(\Lambda)$ }
		
		\For{ $I \subset [d]$}{
			~\\
			Compute the Hermite form $A_I$ of $\Lambda_I$.\\
			$h_I(\Lambda) \leftarrow \val ( \det(A_I) )$ (sum of valuations of diagonal elements of $A_I$)
		}
		$H(\Lambda) \leftarrow (h_I(\Lambda))_{I \subset [d]}$
		
		\caption{Computing $H(\Lambda)$}
		\Return $H(\Lambda)$.
	\end{algorithm}

	Let us now discuss a couple of low-dimensional examples when $K = \Q_p$.
	\begin{example}
		\label{example:2d}
		Let $ \Lambda$ be the lattice represented by $  A = \begin{pmatrix}
			1    & 0        \\
			p    & p^2     \\
		\end{pmatrix} $. The coefficients $h_I(\Lambda)$ of the polynomial $\varphi_\Lambda$ can be computed from the representative $A$ using Algorithm (\ref{algo}) and we have
		\[
		h_\emptyset(\Lambda) = 0, \quad  h_{1}(\Lambda) = 0, \quad h_{2}(\Lambda) = 1, \quad h_{1,2}(\Lambda)= 2	
		\]
		and  then we get
		\[
		\varphi_\Lambda(v_1,v_2) = \max(0, \quad v_1, \quad v_2 - 1, \quad v_1 + v_2 - 2).
		\]
		The independence statement $X_1 \independent X_2$ does not hold since the inequality $h_1(\Lambda) + h_2(\Lambda) \leq h_{12}(\Lambda)$ is strict.
		
		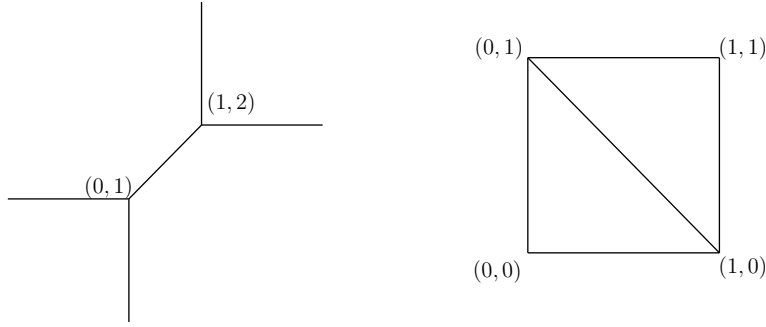
\begin{figure}[H]
		\begin{center}
			\tikzset{every picture/.style={line width=0.5pt}} 
			\begin{tikzpicture}[x=0.75pt,y=0.75pt,yscale=-0.7,xscale=0.7]			
				\draw    (192.65,207) -- (192.65,295.46) ;
				\draw    (192.65,295.46) -- (278.99,295.46) ;
				\draw    (192.65,295.46) -- (140.85,348.54) ;
				\draw    (54.51,348.54) -- (140.85,348.54) ;
				\draw    (140.85,348.54) -- (140.85,437) ;
				\draw    (425.33,387.3) -- (561.95,387.3) ;
				\draw    (425.33,247.07) -- (425.33,387.3) ;
				\draw    (561.95,247.07) -- (561.95,387.3) ;
				\draw    (425.33,247.07) -- (561.95,247.07) ;
				\draw    (561.95,387.3) -- (425.33,247.07) ;
				\draw (107.71,330.94) node [anchor=north west][inner sep=0.75pt]  [xscale=0.7,yscale=0.7]  {$( 0,1)$};			
				\draw (195.12,270.82) node [anchor=north west][inner sep=0.75pt]  [xscale=0.7,yscale=0.7]  {$( 1,2)$};			
				\draw (384.39,391.02) node [anchor=north west][inner sep=0.75pt]  [xslant=-0.05,xscale=0.7,yscale=0.7]  {$( 0,0)$};			
				\draw (386.05,230.76) node [anchor=north west][inner sep=0.75pt]  [xscale=0.7,yscale=0.7]  {$( 0,1)$};			
				\draw (560,388.41) node [anchor=north west][inner sep=0.75pt]  [xscale=0.7,yscale=0.7]  {$( 1,0)$};			
				\draw (560,230.76) node [anchor=north west][inner sep=0.75pt]  [xscale=0.7,yscale=0.7]  {$( 1,1)$};
			\end{tikzpicture}
			\captionof{figure}{ Tropical curve of $\varphi_\Lambda$ and its regular triangulation of the square for example \ref{example:2d}}
		\end{center}
		\end{figure}
	\end{example}
	\begin{example}
		\label{example:3d}
		Let $ \Lambda$ be the lattice represented by $  A = \begin{pmatrix}
			1    & 0      & 0  \\
			1    & \pi^2   & 0  \\
			1    & \pi      & pi^2 
		\end{pmatrix} $. The polynomial $\varphi_\Lambda$ can be computed again using Algorithm (\ref{algo})  and we get 
		\begin{align*}
			&h_\emptyset(\Lambda) = 0   \\
			h_{1}(\Lambda) = 0,  \quad &h_{2}(\Lambda) = 0, \quad  \ \ h_{3}(\Lambda) = 0\\
			h_{1,2}(\Lambda) = 2,  \quad &h_{1,3}(\Lambda) = 1, \quad h_{2,3}(\Lambda) = 1\\
			& h_{1,2,3}(\Lambda) = 4.
		\end{align*}
		So we deduce that  
		\[
		\varphi_{\Lambda}(v) = \max(0, v_1, v_2, v_3,  v_1 + v_2 - 2,  v_1 + v_3 -1 , v_2 + v_3 -1,  v_1 + v_2 + v_3 - 4).
		\]
		We can easily check that the supermodularity inequalities are satisfied. Also, none of the conditional independence statements $X_i \independent X_j | X_k$ are satisfied for $\{i,j,k\}= \{1,2,3\}$ since the point $H(\Lambda)$ is in the interior of the cone $\Scal_3$, i.e, all the inequalities $h_{ki}(\Lambda) + h_{kj}(\Lambda) \leq h_{i}(\Lambda) + h_{ijk}(\Lambda)$ are strict.
		
		\begin{figure}[H]
			\centering
			\begin{subfigure}[b]{0.40\textwidth}
				\includegraphics[scale=0.20]{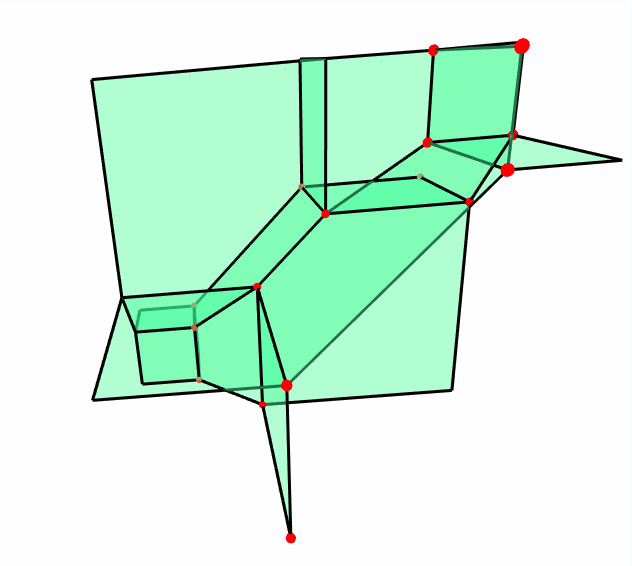}
				\caption{Tropical variety of $\varphi_{\Lambda}$.}
			\end{subfigure}  
			~ 
			\begin{subfigure}[b]{0.40\textwidth}
				\includegraphics[scale=0.40]{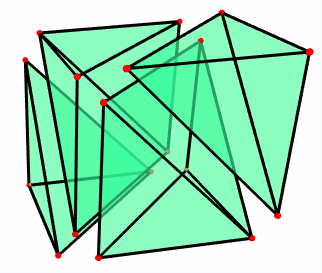}
				\caption{Regular subdivision of the Newton polytope of $\varphi_\Lambda$.}
			\end{subfigure} 
			\caption{ Tropical geometry of the lattice $\Lambda$ for \cref{example:3d}.}
			\label{fig:3d}
		\end{figure}
	\end{example}
	\begin{remark}\label{rem:algo2}
		For any lattice $\Lambda$, there exists a maximal (for inclusion) diagonal lattice inside $\Lambda$ and a minimal diagonal lattice containing $\Lambda$. Let us denote these two lattices by $\bm{\pi^{a}}$ and  $\bm{\pi^{b}}$ respectively, where $a \geq b \in \Z^d$. So, we have the inclusions  $\bm{\pi^{a}} \subset \Lambda \subset  \bm{\pi^{b}}$. It is not difficult to see that the region of linearity corresponding to the monomial $v_1 + \dots + v_d - h(\Lambda)$ in the tropical polynomial $ \varphi_\Lambda(v)$ is the orthant $\R_{\geq a} \coloneqq \{ x \in \R^d, x \geq a \}$. Similarly, the region of linearity corresponding to the monomial $0$ is the orthant $\R_{\le b} \coloneqq \{x \in \R^d, x \leq b\}$. From this, we can the deduce the following  recursive relation
		\[
		h_{[d]}(\Lambda) = h_{[d-1]}(\Lambda) + a_{d}.
		\]
		This iterative way of computing the entropy map $H(\Lambda)$ is slightly more efficient than Algorithm \ref{algo} where we have to compute the whole Hermite form of $\Lambda_I$ for every $I \subset [d]$. This iterative algorithm is the one implemented in (\ref{algo_url}).
	\end{remark}

	\section{The entropy map on non-archimedean fields}\label{Section4}
	
	In this section we generalize some of the results in Section \ref{Section2} to the case where $K$ is a field with a non-archimedean valuation.
	
	When the residue field $k$ of $K$ is infinite or the valuation group $\Gamma$ is dense in $\R$, the probabilistic framework we had in Section \ref{Section2} is no longer valid. More precisely, we lose the local compactness and we no longer necessarily have a Haar measure on $K$. 
	
	We define the entropy map $H$ of a lattice as in Section \ref{Section2}, i.e for any $I \subset [d]$,
	\[
	h_I(\Lambda) \coloneqq \min_{|J| = |I|} \val(\det(A_{I \times J})),
	\]
	where $A$ is a representative of $\Lambda$. We can still define a \emph{Hermite representative } of $\Lambda$.
	\begin{definition}
		
		Every lattice $\Lambda$ in $K^d$ has a representative $A$ in \emph{Hermite normal form}, i.e. a matrix $A$ in $\GL_d(K)$ satisfying the following conditions:
		
		\begin{enumerate}[label=(\roman*)]
			\item $A$ is lower diagonal. 
			\item For any $1\leq j < i\leq d$ we have either $\val(A_{i,j}) < \val(A_{j,j})$ or $A_{i,j} = 0$.
		\end{enumerate}
		
	\end{definition}		 
	
	The same argument used in Theorem \ref{thm:super_mod} can be used again to show that the image of $H$ still lies in the supermodular cone $\Scal_d$. In this setting however, since the valuation group can be dense in $\R$, the image is not necessarily in $\Scal_d \cap \Z^{2^d - 1}$. As in Section \ref{Section2}, the map $H$ fails to be surjective when $d \geq 3$. The algorithm we provide in (\ref{algo_url}) computes the map $H$ when $K = \Q\{\{t\}\}$ is the field of Puiseux series over $\Q$.
	
	Now we show that the only distribution on the field Laurent series $K = \R((t))$ that satisfies the definition suggested in \cite[Definition 4.1]{evans2001local} is the Dirac measure at $0$. Let $\Prob$ be such a probability measure. First, we recall that if $X$ is a random variable with distribution $\Prob$, then for any $a \in \Ocal_{\K}^{\times}$ the random variables $X$ and $a X$ have the same distribution, and we write $X \stackrel{d}{=} a X$. In particular, for any $a \in \R^{\times}$ we have $X \stackrel{d}{=} a X$.
	\begin{proposition}
		The probability distribution $\Prob$ is the Dirac measure at $0$.
	\end{proposition}
	\begin{proof}
		We can write the power series expansion of $X$ as $X = X_0 t^{V} + X_1 t^{V + 1} + \dots$, where $V \in \Z$ is the random valuation of $X$. Hence for $a \in \R^{\times}$ we have $ a X = a X_0 t^{V} + a X_1 t^{V + 1} + \dots $, and we deduce that $X_k \stackrel{d}{=} a X_k$ for any $k \geq 0$ and $a \in \R^{\times}$. We then deduce that $X_k = 0$ almost surely for all $k \geq 0$. Hence $X = 0$ almost surely which finishes the proof.
	\end{proof}
	Using a variant of this argument, it is not difficult to see that a similar problem would arise when we try to define Gaussian measures by orthogonality for all fields listed in Example \ref{ex:fields}. It is not immediately clear how to fix this problem and find a suitable definition for \emph{Gaussian measures} on non-archimedean valued fields.
	
	\begin{Problem}
		Is there a suitable definition for Gaussian measures on the fields listed in Example \ref{ex:fields}?
	\end{Problem}
	
	\begin{remark}
		We can define a probability measure on $\R^d$ induced by $\Lambda$ via its tail distribution $Q_\Lambda $ as in Section \ref{Section2}. One can see that the support of this distribution is $\trop(\Lambda) \coloneqq  \val( \Lambda \cap (K^\times)^d)$; the image under valuation of points in $\Lambda$ with no zero coordinates. This is in general a polyhedral complex in $\R^d$ where each edge is parallel to some $e_I := \sum_{i \in I} e_i$. The following figure is a drawing of  $\trop(\Lambda)$ for a lattice in $K^3$ when $K = \K$ (the field of generalized Puiseux series).
		
		\begin{figure} [H] 
			\begin{center}
				\includegraphics[scale = 0.30]{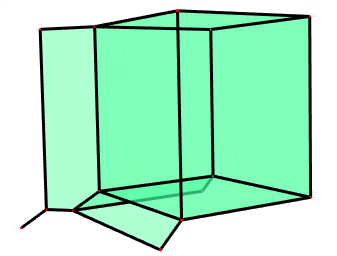}
			\end{center}
				\caption{The polyhedral complex $\trop(\Lambda)$ for $\Lambda$ in \cref{example:3d}.}
				\label{fig:support} 
		\end{figure}
	\end{remark}

	To conclude this section we give a partial answer for Problem \ref{problem} when $d = 2,3$ and the valuation group is $\R$.
	
	\begin{proposition}
		For $d=2$, the image $\im(H) $ of the entropy map $H$ is exactly $\Scal_2$.
	\end{proposition}
	\begin{proof}
		For $\Lambda$ with representative $\begin{pmatrix} t^{a}	 &  0 \\  t^{b}  & t^{b + \delta}\end{pmatrix}$ with $a,b \in \R$ and $\delta \geq 0$ we have $H(\Lambda) = (a,b,a + b + \delta)$. So $H$ is indeed surjective onto $\Scal_2$.
	\end{proof}
	
	For $d=3 $, the cone $\Scal_3 \subset \R^7$ has a lineality space $\Lcal_3$ of dimension $3$. Since both $\Scal_3$ and $\im(H)$ are stable under translations in $\Lcal_3$ (see Remark \ref{rem:lineality} and Lemma \ref{lem: diag_perm} on diagonal scaling of lattices), they are fully determined by their projection onto a complement of $\Lcal_3$. Let us we write vectors $x$ of $\R^7$ in the following form
	\[
	x =(x_1,x_2,x_3 ; \ x_{12}, x_{13}, x_{23}; \ x_{123} ),
	\]	
	and let us project $\Scal_3$ and $\im(H)$ on the linear space $W \subset \R^7$ of vectors of the form 
	\[
	x =(0,x_2,x_3 ; 0, x_{13}, x_{23}; 0).
	\]	
	who is a complement of $\Lcal_3$ in $\R^7$. We write a vector of $W$ as $(x_2,x_3 ;  x_{13} , x_{23}) $ or simply as $(w,x,y,z)$ to simplify notation. Let us denote by $ \Pcal , \Ccal$ be the projections of $\im(H)$ and $\Scal_3$ respectively onto the space $W$. From Section \ref{Section3}, we clearly have $\Pcal \subset \Ccal$.

\begin{figure}
		\centering
		\begin{subfigure}{0.40\textwidth} 
			\includegraphics[width=3cm, height=2.5cm]{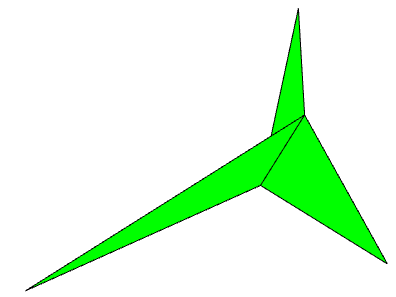}
			\caption{$\Pcal \cap \{ w +x + y +z + 1 = 0\}$.}
		\end{subfigure}  
		~ 
		\begin{subfigure}{0.40\textwidth} 
			\includegraphics[width=4.5cm,height=3cm]{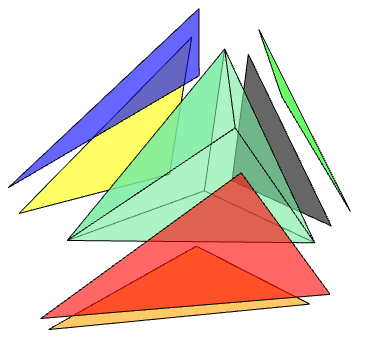}
			\caption{$\Ccal \cap \{  w +x + y +z + 1 = 0\}$. Red facet: $1 \independent 2$; Blue facet: $1 \independent 3$, Green facet: $2 \independent 3$; Orange facet: $1 \independent 2 | 3$; Yellow facet: $1 \independent 3 | 2 $; Grey facet: $2 \independent 3 | 1$.}
		\end{subfigure} 
		\caption{Intersections of $\Pcal$ and $\Ccal$ with the affine hyperplane $x+y+z+w + 1= 0$.}
		\label{fig:Ccal_Pcal}
	\end{figure}

	The projection $\Ccal$ of $\Scal_3$ onto $W$ is a polyhedral cone that does not contains any lines.  In the language of polyhedral geometry, this is called a \emph{pointed cone}. Moreover, the dimension of this projection is $4$. It is defined in $W$ by the inequalities

	\begin{equation}\label{ineq: C}
		\begin{cases}
			w \leq 0, \quad  \quad \quad  x \leq y,  \\ w + x \leq z, \quad 	\ y \leq 0, \\ z \leq w, \quad  \quad  \quad  y + z \leq x.
		\end{cases} 
	\end{equation}
	This defines $\Ccal$ as a pointed cone over a bipyramid (see Figure \ref{fig:Ccal_Pcal}).

	On the other hand, any lattice $\Lambda$ in $\K^3$ can be represented, up to diagonal scaling, by a representative with Hermite form of the shape
	\[
	\begin{pmatrix}
		1    & 0      & 0    \\
		\ast    & 1    & 0    \\
		\ast    & \ast      & 1 
	\end{pmatrix}.
	\]
	The entropy vector of a lattice $\Lambda$ with such a Hermite normal form is of the shape
	\[
	H(\Lambda) = (0,h_2,h_3 ; 0, h_{13}, h_{23}; 0).
	\]
	This corresponds to the projection of $\im(H)$ to $W$ parallel to $\Lcal_3$. So the projection $\Pcal$ of $\im(H)$ onto $W$ is the set
	\[
	\Pcal = \left \{ H(\Lambda), \ \Lambda \text{ given by a matrix of the shape }  \begin{pmatrix}	1    & 0      & 0    \\	\ast    & 1    & 0    \\	\ast    &  \ast      & 1 \end{pmatrix}   \text{ in } \K ^{3 \times 3}  \right\}.
	\]
	For a lattice $\Lambda$ with representative $A = \begin{pmatrix}	1    & 0      & 0    \\	a    & 1    & 0    \\	b    &  c      & 1 \end{pmatrix}$, such that $a,b,c \in  \K$ with negative or zero valuation (see Definition \ref{Hermite}), the point $H(\Lambda)$ in $W$ is given by
	\[
	\begin{cases}
		w = h_2(\Lambda) = \val(a) , \\  x =  h_3(\Lambda) = \min( \val(b), \val(c) ), \\  y = h_{13}(\Lambda) = \val(c),  \\ z=  h_ {23}(\Lambda) = \min(\val(ac -  b ) , \val(a)).
	\end{cases}
	\]
	One can check that, for any choice of $a,b,c \in \K$ with negative or zero valuation, the above coordinates satisfy the inequalities in (\ref{ineq: C}). With the constraints on the valuations of $a,b,c$, and from this parametric representation of $\Pcal$, we can see that points of $\Pcal$ have to satisfy the inequalities
	\[
	\begin{cases}
		w \leq 0, \\ x \leq y, \\ y \leq 0.
	\end{cases}
	\]
	The only part that remains to determine is the inequalities involving the last variable $z$. The ambiguity comes from the fact that cancellations can happen in $ac - b$ which might affect $\val(ac - b)$ and hence also $z$. But, separating the cases where $\val(ac) = \val(b)$ and $\val(ac) \neq \val(b)$, we get the following three sets of inequalities that describe $\Pcal$ as a polyhedral complex:
	\[
	\begin{cases}
		w \leq 0, \\ x \leq w + y, \\ y \leq 0, \\ z = x,
	\end{cases} \text{ ,  } \quad  
	\begin{cases}
		w \leq 0, \\ x \leq y, \\ y \leq 0, \\ y+w \leq x, \\  z = y+w,
	\end{cases} \text{ and } \quad 
	\begin{cases}
		w \leq 0, \\ y \leq 0, \\  x = y + w,  \\ z \leq w, \\  x  \leq z.
	\end{cases}
	\]

	We can then see that $\Pcal$ is a polyhedral fan of dimension $3$ inside $\Ccal$. More precisely, $\Pcal$ is the union of  three  pointed polyhedral cones of dimension $3$ inside $\Ccal$ which is a cone of dimension $4$. \cref{fig:Ccal_Pcal} depicts the intersections of $\Pcal$ and $\Ccal$ with the hyperplane $w + x + y + z + 1 = 0$ (slicing the pointed cones with a hyperplane).  
	
	\begin{corollary}
		The entropy map $H : \GL_d(\K) / \GL_d(\Ocal_{\K}) \to \Scal_d$  is not surjective when $d \geq 3$.
	\end{corollary}

	We expect this result to hold in every dimension, i.e, the image $\im(H)$ is a polyhedral fan whose facets are polyhedral cones of dimension $\frac{d(d+1)}{2}$ inside $\Scal_d$ which is of dimension $2^{d} - 1$.

	\section{Conclusion}
	In conclusion, there are many similarities between the classical theory of Gaussian distributions on euclidean spaces and the theory of Gaussian measures on local fields as defined by Evans in \cite{evans2001local}. In this paper we have exhibited another similarity in terms of differential entropy. This gives reason to think that the suitable non-archimediean analog of the positive definite cone is indeed the set of lattices, or more precisely, in the language of group theorists, the Bruhat-Tits building for $\SL$. This analogy can still be carried out for non-archimedean valued fields in general. However, when the field $K$ has a dense valuation group or an infinite residue field, we lose the probabilistic interpretation and thus also the notion of entropy.

		\bibliographystyle{alpha}
		\bibliography{article.bib}
		
	\end{document}